\documentclass[10pt]{amsart}
\usepackage{amsfonts}
\usepackage{mathrsfs}
\usepackage{amssymb}
\usepackage{amsmath}
\usepackage{epic}

\newtheorem{proposition}{Proposition}[subsection]
\newtheorem{lemma}[proposition]{Lemma}
\newtheorem{corollary}[proposition]{Corollary}
\newtheorem{theorem}[proposition]{Theorem}

\theoremstyle{definition}

\theoremstyle{remark}
\newtheorem{remark}[proposition]{Remark}

\def \lhp{\hspace*{-1mm}\leftharpoonup\hspace*{-1mm}}

\newcommand{\thlabel}[1]{\label{th:#1}}
\newcommand{\thref}[1]{Theorem~\ref{th:#1}}
\newcommand{\selabel}[1]{\label{se:#1}}
\newcommand{\seref}[1]{Section~\ref{se:#1}}
\newcommand{\lelabel}[1]{\label{le:#1}}
\newcommand{\leref}[1]{Lemma~\ref{le:#1}}
\newcommand{\prlabel}[1]{\label{pr:#1}}
\newcommand{\prref}[1]{Proposition~\ref{pr:#1}}
\newcommand{\colabel}[1]{\label{co:#1}}

\newcommand{\relabel}[1]{\label{re:#1}}
\newcommand{\reref}[1]{Remark~\ref{re:#1}}

\newcommand{\eqlabel}[1]{\label{eq:#1}}
\newcommand{\equref}[1]{(\ref{eq:#1})}

\def\cd{\cdot}
\def\cds{\cdots}

\def\hp{\hspace{1.5cm}}
\def\lhp{\hspace{3.5cm}}

\def\<{\leq}
\def\>{\geq}
\def\a{\alpha}
\def\b{\beta}
\def\d{\delta}

\def\D{\Delta}
\def\g{\gamma}

\def\ep{\varepsilon}
\def\om{\omega}

\def\vs{v_{\s}}

\def\k{\textbf{k}}
\def\l{\lambda}

\def\Ls{L(\s)}

\def\ms{M(\s)}
\def\s{\sigma}

\def\ti{\times}
\def\op{\oplus}
\def\ot{\otimes}
\def\vp{\varpi}

\def\ra{\rightarrow}

\def\mcg{\mathcal{G}}
\def\mcw{\mathcal{W}}

\def\mfg{\mathfrak{g}}
\def\Z{\mathbb{Z}}

\def\Hom{\mathrm{Hom}}

\def\qg{\mathcal{U}}
\def\qgg{\mathcal{U}^{\>0}}

\def\yd{\mathcal {Y}\mathcal {D}}

\date{}

\begin{document}
\title{The Half-quantum Group $\qgg$}
\thanks{The work is supported by NSF, No.10471121}

\author{Zhen Wang}
\address{Department of Mathematics, Yangzhou University,
Yangzhou 225002, China} \email{wangzhen118@gmail.com}
\author{Hui-Xiang Chen}
\address{Department of Mathematics, Yangzhou University,
Yangzhou 225002, China} \email{hxchen@yzu.edu.cn} 

\begin{abstract}
Let $\qgg$ denote the half quantum group for a fixed simple Lie
algebra. We examine some properties and representation of $\qgg$. We
prove that the Hopf algebra $\qgg$ is not quasi-cocommutative, and
hence the category of left $\qgg$-module is not a braided monoidal
category. In the weight module category $\mcw$, we describe all the
simple objects and the projective objects. We also describe all
simple Yetter-Drinfel'd  $\qgg$-weight modules.
\end{abstract}
\maketitle

\section*{\bf Introduction}
For a semisimple Lie algebra $\mathfrak{g}$, Drinfel'd and Jimbo
introduced a class of non-commutative non-cocommutative Hopf
algebras $U_q (\mathfrak{g})$ independently(see \cite{dri,jim}).
These algebras are called the quantized enveloping algebra and
regarded as certain deformations of the enveloping algebra
$U(\mathfrak{g})$. There are many authors to have studied them. The
reader can find a detailed introduction about them in any one of
\cite{dck,jan,kas}. Let $\qgg$ denote the upper triangular Hopf
subalgebra of $U_q (\mathfrak{g})$, called the half quantum group.
Then the Yetter-Drinfel'd modules over $\qgg$ provide universal
$R$-matrices and solutions to the quantum Yang-Baxter equation. They
are also of interest in connection with knot theory and invariant of
three-manifolds(\cite{maj,ret}).

When $q$ is not a root of unity, the representation theory of $U_q
(\mathfrak{g})$ is essentially the same as that of
$U(\mathfrak{g})$. However, if $q$ is a root of unity, the situation
changes dramatically since $U_q (\mathfrak{g})$ contains a large
central subalgebra in this case(\cite{dck}). For the half quantum
group $\qgg$, the similar situation appears. In this case, one can
form a finite dimensional quotient Hopf algebra $u^{\>0}$ of $\qgg$.
In \cite{gun}, Gunnlaugsd\'ottir discussed the monoidal property of
$u^{\>0}$ when $\mathfrak{g}=\mathfrak{sl}_2$. For simply laced Lie
algebra $\mathfrak{g}$, Cibils gave a presentation of $u^{\>0}$ by
quiver and relations, and showed that only
$u^{\>0}(\mathfrak{sl}_2)$ is of finite representation type, the
others are of wild representation type \cite{cib}.

The purpose of this article is to examine some properties and
representation theory of $\qgg$ for any simple Lie algebra
$\mathfrak{g}$.

In section $1$, we first review the definition and some properties
of quantum group $\qg=U_q(\mathfrak{g})$, and give the definition
and some properties of $\qgg$. Then we discuss the
quasi-cocommutative property of a graded Hopf algebra. We show that
$\qgg$ is not quasi-cocommutative, and hence the category
$_{\qgg}\mathcal {M}$ of left $\qgg$-modules is not a braided
monoidal category. Recall that the category $_H\yd ^H$ of
Yetter-Drinfel'd $H$-modules is a braided monoidal category when
$H^{cop}$ is a Hopf algebra (cf. \cite{yet}). It is well-known that
the center ${\mathcal Z}(_H{\mathcal M})$ of $_H{\mathcal M}$ is a
braided monoidal category, tensor equivalent to $_H\yd ^H$ (cf.
\cite[Theorem XIII.5.1]{kas}). We also consider the quotient algebra
$u^{\>0}(\mathfrak{g})$ and generalize the results in \cite{gun}. In
section $2$, we discuss the simple modules, Verma modules and the
indecomposable projective objects in the weight module category
$\mcw$. Moreover, the monoidal structure of $\mcw$ is considered. We
give the decomposition of the tensor product of two Verma modules.
We also describe the simple Yetter-Drinfel'd $\qgg$-modules which
are weight modules. Using Radford's results in \cite{rad}, we show
that there is a one-to-one correspondence between the set $\mcg
((\qgg)^\circ)\times \mcg (\qgg)$ and the set $\mathcal{E}$ of
isomorphic classes of the simple Yetter-Drinfel'd $\qgg$-modules
which are weight modules.

For basic background about quantum group and Hopf algebra, the
reader is directed to \cite{jan,kas,mon}.

\section{\bf The structure of $\qgg$}\selabel{1}
\subsection{}\selabel{1.1}
Throughout, $k$ is an algebraically closed field with characteristic
$0$, $k^{\times}=k\setminus \{0\}$, and $q\in k^{\times}$ with
$q\neq\pm 1$. Unless otherwise stated, all algebras, Hopf algebras
and modules are defined over $k$; dim, $\otimes$ and Hom stand for
${\rm dim}_k$, $\otimes_k$ and ${\rm Hom}_k$, respectively. We will
use Sweedler's sigma notation for the comultiplicatin  of a Hopf
algebra (cf. \cite{mon}). Let $\Z$ denote the integer set, and
$\Z_+$ denote the non-negative integer set. For $n\in \Z$, let
$[n]_q=(q^{n}-q^{-n})/(q-q^{-1})$. As usual, we define $[0]_q!=1$
and $[n]_q!=[n]_q[n-1]_q\cdots \cdots [1]_q$ for $n\geq 1$, and the
Gaussian $q$-binomial coefficients
$$\left [ \begin{array}{c}n\\ j \end{array} \right ]_{q}
=\frac{[n]_q!}{[j]_q![n-j]_q!}\ , \hspace{1cm} n\geq j\geq 0.$$ For
a fixed simple Lie algebra $\mathfrak{g}$ with rank $n$, let
$C=(a_{ij})_{n\ti n}$ be its Cartan matrix. Then there exists a
diagonal matrix $D=diag\{d_1,d_2,\cdots,d_n\}$ over $\Z$ such that
$DC$ is symmetric, i.e., $d_ia_{ij}=d_ja_{ji}$ for $1\<i,j\<n$. We
may assume that each $d_i>0$ and $\sum _i d_i$ is as minimal as
possible. Let $q_i=q^{d_i}$. Then the quantized enveloping algebra
$\qg=U_q (\mathfrak{g})$ associated to $\mathfrak{g}$ is a
$k$-algebra with generators $E_i, F_i , K_i, K_i^{-1}(1\< i\<n)$
subject to the relations: for all $1\<i,j\<n$,
\begin{eqnarray}
&&K_iK_j=K_jK_i\ ,\, K_iK_i^{-1}=K_i^{-1}K_i=1\ , \\
&&K_iE_jK_i^{-1}=q_j^{a_{ji}}E_j\ , \\
&&K_iF_jK_i^{-1}=q_j^{-a_{ji}}F_j\ ,\\
&&E_iF_j-F_jE_i=\d _{ij}(K_i-K_i^{-1})/(q_i-q_i^{-1})\ ,\\
&&\sum_{s=0}^{1-a_{ij}}(-1)^s\left [ \begin{array}{c}1-a_{ij}\\ s
\end{array} \right ]_{q_i}E_i^{1-a_{ij}-s}E_jE_i^s=0\ , \quad \text{if
}i\neq j\ ,\\
&&\sum_{s=0}^{1-a_{ij}}(-1)^s\left [ \begin{array}{c}1-a_{ij}\\ s
\end{array} \right ]_{q_i}F_i^{1-a_{ij}-s}F_jF_i^s=0\ , \quad \text{if }i\neq j\ .
\end{eqnarray}

$\qg $ is a Hopf algebra with comultiplication $\D$, antipode $S$
and counit $\ep$ defined by
\begin{eqnarray}
&\D E_i=E_i\ot 1+K_i\ot E_i\ ,\,\D F_i=F_i\ot K_i^{-1}+1\ot F_i\ ,
    \, \D K_i=K_i\ot K_i\ ,\\
&SE_i=-K_i^{-1}E_i\ ,\,SF_i=-F_iK_i\ ,\,SK_i=K_i^{-1}\ ,\\
&\ep E_i=0\ ,\,\ep F_i=0\ ,\, \ep K_i=1\ ,
\end{eqnarray}
where $i=1, 2, \cdots, n$.

Let $\qg ^+$, $\qg^-$ and $\qg ^0$ be the subalgebras of $\qg$
generated by $\{E_i|1\leq i\leq n\}$, $\{F_i|1\leq i\leq n\}$ and
$\{K_i,K_i^{-1}|1\< i\< n\}$ respectively. It follows from (1)-(6)
that $\qg =\qg ^-\qg ^0\qg ^+$. Moreover, the multiplication gives a
$k$-vector space isomorphism
\begin{equation}
\qg ^-\ot \qg ^0\ot \qg ^+\cong \qg\ .
\end{equation}

Let $\qgg=\qg ^0 \qg ^+$ be a subalgebra of $\qg$ generated by
$\{E_i,\,K_i,K_i^{-1}|1\< i\< n\}$. Actually, we can directly define
$\qgg$ as an $k$-algebra with the generators $E_i,\,K_i,K_i^{-1},
1\< i\< n$, and the relations (1),(2) and (5) above. Let $(\qgg)^+$
and $(\qgg)^0$ denote the subalgebras of $\qgg$ generated by
$\{E_i|1\<i\<n\}$ and $\{K_i,K_i^{-1}|1\< i\< n\}$, respectively.
Then clearly $(\qgg)^+=\qg ^+$ and $(\qgg)^0=\qg ^0$. It's easy to
check that $\qgg$ is a Hopf subalgebra of $\qg$. In this paper, we
will mainly study the properties and the representation theory of
$\qgg$.

\subsection{}\selabel{1.2}
For a given Cartan matrix $(a_{ij})$, let $P=\sum_{i=1}^n\Z \vp_i$
be the weight lattice. Define simple roots by
$$ \a_j=\sum_{i=1}^n a_{ij}\vp_i, \quad
j=1,\cdots,n.$$

Let $\D=\{\a_1,\cdots ,\a_n\}$, $Q=\Z\D$ (the root lattice), and
$Q_+=\sum_i \mathbb{Z}_+\a_i$. Then there is a partial ordering on
$P$ defined by $\mu\<\l$ if $\l -\mu \in Q_+$.

Define automorphisms $\g _i$ of $P$ by $\g _i\vp _j=\vp _j-\d_{ij}\a
_i\, (i,j=1,\cdots,n).$ Then $\g_i \a_j=\a_j-a_{ij}\a_i$. Let $W$ be
the (finite) subgroup of $GL(P)$ generated by $\g_1,\cdots, \g_n$,
called the Weyl group. Then $Q$ is $W$-invariant. Let
$R=W\D,\,R^+=R\cap Q_+$ and $R^-=-R^+$. Then $R$ is a root system
corresponding to the Cartan matrix $(a_{ij})$, $R^+$ a set of
positive roots, $R=R^+\cup R^-$.

Fix a reduced expression $\g_{i_1}\g_{i_2}\cdots \g_{i_N}$ of the
longest element $\om _0$ of $W$. This gives us an ordering of the
set of positive roots $R^+$:
$$\b _1=\a_{i_1},\, \b_2=\g_{i_1}\a_{i_2},\,\cdots,\,
\b_N=\g_{i_1}\cdots \g_{i_{N-1}}\a_{i_N}.$$

Let $T_i\, (i=1,2,\cds,n)$ be the automorphisms of $\qg$ satisfying
$$T_i(E_j)=\left\{\begin{array}{ll}
-F_iK_i,&\text{if
}i=j,\\
\sum_{l=0}^r (-1)^{l+r}q_i^{-l}E_i^{(r-l)}E_jE_i^{(l)},& \text{if
}i\neq j,\\
\end{array}\right.$$
where $r=-a_{ij}$ and $E_i^{(l)}=\frac{E_i^l}{[l]_{q_i}!}$. Define
root vectors(see \cite{dck,jan}) by
$$E_{\b_s}=T_{i_1}\cdots T_{i_{s-1}}E_{i_s},\quad 1\<s\<N.$$
From now on, let $\mathfrak{g}$ be a simple Lie algebra with rank
$n$, $C=(a_{ij})_{n\ti n}$ be its Cartan matrix, and
$D=diag\{d_1,d_2,\cdots,d_n\}$ be a diagonal matrix over $\Z$ such
that $DC$ is symmetric as before. We also assume that each $d_i>0$
and $\sum_{i=1}^nd_i$ is minimal.

\begin{theorem}\thlabel{1.2.1} With
the above notations, we have\\
$(1)$ $\{E_{\b _1},E_{\b _2},\cdots ,E_{\b _N}\}\subset (\qgg)^+$.\\
$(2)$ $\{E_1,E_2,\cdots ,E_n\}\subset \{E_{\b _1},E_{\b _2},\cdots
       ,E_{\b _N}\}$. More precisely, suppose that $\a _t =\b _j$.
\mbox{\hspace{0.5cm}}Then $E_{\b _j}=E_t$.\\
$(3)$ For an integer $r$ with $1\<r\<N$, let
$\om=\g_{i_1}\g_{i_2}\cdots\g_{i_r}\in W$ with $l(\om)=r$. Let\\
\mbox{\hspace{0.5cm}}$(\qgg)^+(\om)=k\langle E_{\b _1},E_{\b
_2},\cdots,E_{\b _r}\rangle$ be the subalgebra of $(\qgg)^+$
generated by\\
\mbox{\hspace{0.5cm}}$E_{\b _1},E_{\b _2},\cdots,E_{\b _r}$. Then
the monomials $E_{\b _1}^{m_1}E_{\b _2}^{m_2}\cdots E_{\b _r}^{m_r}$
are a basis of $(\qgg)^+(\om)$, \mbox{\hspace{0.5cm}}where $m_1,m_2,\cdots,m_r$ run over $\Z _+$.\\
$(4)$ $(\qgg)^+(\om _0)=(\qgg)^+$.\\
$(5)$ $\qgg\simeq (\qgg)^0\ot (\qgg)^+$ as vector spaces.
\end{theorem}
\begin{proof}
It follows from the properties of $\qg$ given in \cite{brg} and the
fact $(\qgg)^+=\qg ^+, \, (\qgg)^0 =\qg ^0$.
\end{proof}

We may also describe $\qgg$ as follows. $(\qgg)^+$ is a right
$(\qgg)^0$-module algebra with the right action of $(\qgg)^0$ on
$(\qgg)^+$ determined by
$$ E_j\cd K_i=q_j^{-a_{ji}}E_j,\ 1\<i,j\<n.$$
Then one can construct the smash product algebra $(\qgg)^0
$\#$(\qgg)^+ $: $(\qgg)^0$\#$(\qgg)^+ =(\qgg)^0\ot (\qgg)^+$ as a
vector space, and the multiplication is given by
$$ (a\ot x)(b\ot y)=\sum ab_{(1)}\ot (x\cd b_{(2)})y,\quad a,b \in
(\qgg)^0,\ x,y\in(\qgg)^+.$$

\begin{theorem}\thlabel{}
$\qgg\cong (\qgg)^0$\# $(\qgg)^+$ as algebras.
\end{theorem}
\begin{proof}
It is a straightforward verification.
\end{proof}

For $\textbf{k}=(k_1,\,\cdots,\,k_N)\in \Z_+^N$ and $\l
=\sum_{i=1}^nt_i\a_i \in Q$, set $E^\k=E_{\b_1}^{k_1}\cdots
E_{\b_N}^{k_N} \text{ and }K_{\l}=K_1^{t_1}\cdots K_n^{t_n}$. Then
by \thref{1.2.1}, $\qgg$ has a PBW basis given by
$\{E^{\k}K_{\l}|\k\in \Z_+^N,\ \l \in Q\}$ .

For $1\<i\<n$, let $deg(K_i)=0$ and $deg(E_i)=\a_i$. Then $\qgg$ is
a $Q_+(=\Z_+ \D)$-graded algebra since the definition relations of
$\qgg$ are homogeneous polynomials under the grading. Obviously,
$\qgg$ is a graded Hopf algebra with respect to the grading.

Let $I$ be a fixed set. Then we have the following lemma.
\begin{lemma}\lelabel{1}
Let $A=\bigoplus _{i\in \Z_+I}A_i$ be a $\Z_+I$-graded algebra. If
$x=\sum_{i\in \Z_+I}x_i \in A$ is an invertible element, then so is
$x_0 \in A_0$. Consequently, $x_0$ is nonzero.
\end{lemma}
\begin{proof}
Let $y=\sum y_i \in A$ such that $xy=yx=1$. Then
$xy=(x_0+x')(y_0+y')=x_0y_0+x_0y'+x'y_0+x'y'=1$, where
$x'=\sum_{i\neq 0}x_i$ and $y'=\sum_{i\neq 0}y_i$. By comparing the
homogeneous components of degree 0 of the two sides, one gets
$x_0y_0=1$. Similarly, $y_0x_0=1$. Therefore, $x_0$ is invertible.
\end{proof}
Recall that a Hopf algebra $H$ is called quasi-cocommutative if
there exists an invertible element $R \in H\ot H$ such that $\D
^{op}(h)=R\D (h)R^{-1}$ for any $h\in H$. In this case, we also say
that $(H, R)$ is quasi-cocommutative.

\begin{lemma}\lelabel{2}
Let $H$ be a $\Z _+I$-graded Hopf algebra. Then $H\ot H$ is also a
$\Z _+I$-graded Hopf algebra with $(H\ot
H)_i=\bigoplus_{j+l=i}H_j\ot H_l$ for all $i\in\Z_+I$. Assume there
exists an invertible element $R\in H\ot H$ such that $(H,R)$ is
quasi-cocommutative. Let $R=\sum_i R_i$ with $R_i\in(H\ot H)_i$,
$i\in\Z_+I$. Then both $(H,\,R_0)$ and $(H_0,\,R_0)$ are
quasi-cocommutative.
\end{lemma}
\begin{proof} Note that $(H\ot H)_0=H_0\ot H_0$. It follows from
\leref{1} that $R_0$ is invertible. Let $h\in H_i$ with $i\in\Z_+I$.
Since $(H,R)$ is quasi-cocommutative, we have
$\D^{op}(h)=R\D(h)R^{-1}$, and hence $\D^{op}(h)R=R\D(h)$. By
comparing the homogeneous components of degree $i$ of the two sides
of the last equation, one gets $R_0\D(h)=\D^{op}(h)R_0$. Hence
$(H,R_0)$ is quasi-cocommutative. Since $H_0$ is Hopf subalgebra of
$H$ and $R_0$ is an invertible element in $H_0$, $(H_0, R_0)$ is
also quasi-cocommutative.
\end{proof}

Using the above lemma, one can get the following theorem.

\begin{theorem}\thlabel{1}
The Hopf algebra $\qgg$ is not quasi-cocommutative.
\end{theorem}
\begin{proof}
If there exists an invertible $R\in \qgg \ot \qgg$ such that
$(\qgg,\, R)$ is quasi-cocommutative, then by \leref{2} we may
assume $R\in (\qgg)_0\ot (\qgg)_0=(\qgg)^0 \ot (\qgg)^0 $ . Let
$R=\sum_{j,l\in\Z ^n}a_{j,l}K^j\ot K^l$, where
$K^r=K_1^{r_1}K_2^{r_2}\cds K_n^{r_n}$ for $r=(r_1,r_2,\cds,r_n)\in
\Z ^n$. Let $e_i=(\d _{1i},\,\cdots,\,\d_{ni}) \text{ for }
i=1,\,\cdots,\, n$. Putting $h=E_1,\, \cdots,\,E_n$ in the equation
$R\D (h)=\D ^{op}(h)R$, respectively, one gets
\begin{equation} \eqlabel{eq:1}
\begin{aligned}
        &\sum_{j,l\in\Z^n} a_{j,l}q_1^{\sum_{i=1}^n a_{1i}j_i}E_1K^j\ot K^l
            +\sum_{j,l\in\Z^n} a_{j,l}q_1^{\sum_{i=1}^n a_{1i}l_i}K^{j+e_1}\ot E_1K^l\\
            &\lhp=\sum_{j,l\in\Z^n} a_{j,l}K^j\ot
                      E_1K^l+\sum_{j,l\in\Z^n} a_{j,l}E_1K^j \ot K^{l+e_1} ,\\
&\sum_{j,l\in\Z^n} a_{j,l}q_2^{\sum_{i=1}^n a_{2i}j_i}E_2K^j\ot K^l
            +\sum_{j,l\in\Z^n} a_{j,l}q_2^{\sum_{i=1}^n a_{2i}l_i}K^{j+e_2}\ot E_2K^l\\
            &\lhp=\sum_{j,l\in\Z^n} a_{j,l}K^j\ot
                      E_2K^l+\sum_{j,l\in\Z^n} a_{j,l}E_2K^j \ot K^{l+e_2} ,\\
   &\lhp \cdots \hp \cdots\\
   &\sum_{j,l\in\Z^n} a_{j,l}q_{n-1}^{\sum_{i=1}^n a_{n-1,i}j_i}E_{n-1}K^j\ot K^l
            +\sum_{j,l\in\Z^n} a_{j,l}q_{n-1}^{\sum_{i=1}^n a_{n-1,i}l_i}K^{j+e_{n-1}}\ot E_{n-1}K^l\\
            &\lhp=\sum_{j,l\in\Z^n} a_{j,l}K^j\ot
                      E_{n-1}K^l+\sum_{j,l\in\Z^n} a_{j,l}E_{n-1}K^j \ot K^{l+e_{n-1}} ,\\
   &\sum_{j,l\in\Z^n} a_{j,l}q_n^{\sum_{i=1}^n a_{ni}j_i}E_nK^j\ot K^l
            +\sum_{j,l\in\Z^n} a_{j,l}q_n^{\sum_{i=1}^n a_{ni}l_i}K^{j+e_n}\ot E_nK^l\\
            &\lhp=\sum_{j,l\in\Z^n} a_{j,l}K^j\ot
                      E_nK^l+\sum_{j,l\in\Z^n} a_{j,l}E_nK^j \ot K^{l+e_n} .\\
                          \end{aligned}
                          \end{equation}
Comparing the coefficients of these equations, we have
\begin{equation} \eqlabel{eq:2}
\begin{aligned}
         a_{j,l} &= q_1^{\sum_{i=1}^n a_{1i}l_i}a_{j-e_1,l},\hp
                     a_{j,l} = q_1^{\sum_{i=1}^n a_{1i}j_i}a_{j,l+e_1} ,\\
         a_{j,l} &= q_2^{\sum_{i=1}^n a_{2i}l_i}a_{j-e_2,l},\hp
                     a_{j,l} =q_2^{\sum_{i=1}^n a_{2i}j_i}a_{j,l+e_2} ,\\
      &\cdots \hp \cdots\hp \cdots\hp\cdots \\
         a_{j,l} &= q_{n-1}^{\sum_{i=1}^n a_{n-1,i}l_i}a_{j-e_{n-1},l},\quad
                     a_{j,l} =q_{n-1}^{\sum_{i=1}^n a_{n-1,i}j_i}a_{j,l+e_{n-1}} ,\\
         a_{j,l} &= q_n^{\sum_{i=1}^n a_{ni}l_i}a_{j-e_n,l},\hp
                     a_{j,l} = q_n^{\sum_{i=1}^n a_{ni}j_i}a_{j,l+e_n} .\\
                      \end{aligned}
                          \end{equation}
Therefore, for any $j, l, j', l'\in\Z^n$, we have $a_{j,
l}=ba_{j',l'}$ for some $b\in k^{\times}$. Since $R=\sum_{j,l\in\Z
^n}a_{j,l}K^j\ot K^l$ is a finite sum (i.e., almost all $a_{j,l}$
are equal to zero), all $a_{j,l}$ are equal to zero. So $R=0$, which
contradicts to the invertibility of $R$.

\end{proof}
When $q$ is a primitive $r$-th root of unity, set
$$d=\left\{\begin{array}{ll}
    r, & \text{if } r \text{ is odd}, \\
    r/2,& \text{if }r \text{ is even}.
  \end{array}\right.$$
Let $u^{\>0}$ be the quotient algebra of $\qgg$ modulo the ideal
generated by $\{K_i^d-1,E_i^d|1\<i\<n\}$. If
$d>d_0=\text{max}\{d_1,d_2,\cds,d_n\}$, then $u^{\>0}$ is a graded
Hopf algebra, where $d_1,d_2,\cds,d_n$ are the diagonal entries of
$D$. If $d\leq d_0$, then $u^{\>0}$ is not a Hopf algebra, i.e., the
ideal generated by $\{K_i^d-1,E_i^d|1\<i\<n\}$ is not a Hopf ideal.
In the rest of this section, we always assume $d>d_0$.

For $j=(j_1,j_2,\cdots,j_n)$, $l=(l_1,l_2,\cdots,l_n)\in(\Z/d\Z)^n$
and $X=(x_{si})_{n\times n}\in M_n(\Z)$, let $A(j,l)=jXl^t=\sum
_{s,i=1}^n x_{si}j_sl_i\in \Z/d\Z$, where $l^t$ is the transpose of
$l$. Clearly, $A(j,l)=A(l,j)$ if $X$ is a symmetric matrix.

\begin{theorem}\thlabel{th:1}
$u^{\>0}$ is quasi-cocommutative if and only if
$\mfg=\mathfrak{sl}_2\text{ and } r= 4.$
 \end{theorem}

\begin{proof}
If $\mfg=\mathfrak{sl}_2\text{ and } r= 4$, then
$u^{\>0}(\mathfrak{sl}_2)$ is exactly the Sweedler's $4$-dimensional
Hopf algebra. It's well-known that $(u^{\>0}(\mathfrak{sl}_2),R)$ is
a triangular Hopf algebra, where $R=1/2(1\ot 1+1\ot K+K\ot 1-K\ot
K)$. Hence $(u^{\>0}(\mathfrak{sl}_2),R)$ is quasi-cocommutative.

Now we assume $\mfg\neq \mathfrak{sl}_2$ or $r\neq 4$, and assume
there exists an invertible element $R=\sum_{j,l\in( \Z/d\Z)
^n}a_{j,l}K^j\ot K^l \in (u^{\>0})_0\ot (u^{\>0})_0$ such that
$(u^{\>0},R)$ is quasi-cocommutative. An argument similar to the
proof of \thref{1} shows that we have equations \equref{eq:2}, where
$j,l\in(\Z/d\Z)^n$. It follows from \equref{eq:2} that
$a_{j,0}=a_{0,l}=a_{0,0}$ and
$a_{j,0}=q^{A(j,l)}a_{j,l}=q^{2A(j,l)}a_{0,l}$ for all
$j,l\in(\Z/d\Z)^n$, where $A(j,l)=jDCl^t\in\Z/d\Z$ as above, and
$DC$ is the symmetrization of the corresponding Cartan matrix $C$.
Now one can deduce a contradiction for each case by considering the
value of $A(j,l)$. For example, if $\mfg $ is of type $G_2$, then
$n=2$ and $DC=\begin{pmatrix}6&-3\\-3&2\end{pmatrix}$. In this case,
let $j=(0,1),l=(1,1)$. Then $A(j,l)=-1$. Since $q^{-2}\neq 1$,
$a_{0,l}=a_{j,0}=a_{0,0}=0$, and so $a_{j,l}=0$ for all
$j,l\in(\Z/d\Z)^2$. Thus $R=0$, a contradiction.
\end{proof}

For the case of $\mfg=\mathfrak{sl}_2\text{ and }r= 4$, it follows
from the proof of \thref{th:1} that an invertible element $R\in
(u^{\>0})_0\ot (u^{\>0})_0 $ with $(u^{\>0}(\mfg),R)$ being
quasi-cocommutative must have the form
$$R=\sum_{j,l\in \Z/d\Z}q^{A(j,l)}aK^j\ot K^l=\sum_{j,l\in
\Z/d\Z}q^{2jl}aK^j\ot K^l,\, \forall a\in k^{\times}.$$

\begin{remark}
1) For the case $u^{\>0}(\mathfrak{sl}_2)$, Gunnlaugsd\'ottir proved in \cite{gun}
that $u^{\>0}(\mathfrak{sl}_2)$ is quasi-cocommutative if and only if $r=4$.\\
2) In \cite{kas}, it was proved that the quantized enveloping
algebra $U_h(\mfg)$ is quasi-cocommutative, and so is the quotient
Hopf algebra of $U_q(\mfg)$ modulo the ideal generated by $\{K_i^d
-1,E_i^d,F_i^d|1\<i\<n \}$ (which is a Hopf ideal of $U_q(\mfg)$)
when $q$ is a root of unity of order $r$.
\end{remark}

\section{\bf The representation of $\qgg$}\selabel{2}
In this section, we discuss the representation theory of $\qgg$. We
use the notations in \seref{1} and assume that $q$ is not a root of
unity.

Let $H$ be a Hopf algebra with antipode $S$, and $H^*=\Hom_k(H,k)$
be the dual algebra of $H$. Let $H^\circ$ denote the finite dual of
$H$, i.e., $H^\circ =\{f\in H^*| f(I)=0\text{ for some ideal } I
\text{ with }{\rm dim}(H/I)<\infty\}$. Then $H^\circ$ has an induced
Hopf algebra structure. Let $M$ and $N$ be two left modules over
$H$. Then $M^*=\Hom_k(M,k)$ and $M\ot _kN$ are also left $H$-modules
defined by $(h\cd f)(m)=f(S(h)\cd m)$ and $h\cd(m\ot n)=h_{(1)} \cd
m\ot h_{(2)} \cd
 n$, respectively, where $h\in H$, $f\in M^*$, $m\in M$, $n\in N$ (see \cite{mon}).

\subsection{The simple modules}\selabel{}

Let $_H\mathcal {M}$ denote the category of all the left modules
over a Hopf algebra $H$. In what follows, an $H$-module means a left
$H$-module, and an $H$-comodule means a right $H$-comodule.

Let $M$ be a $\qgg$-module. For any
$\s=(\s_1,\s_2,\cdots,\s_n)\in(k^{\times})^n$, let $M_{(\s)}=\{v\in
M| K_i\cdot v=\s_iv, 1\leq i\leq n\}$. An element $v\in M$ is called
a weight vector if $v\in M_{(\s)}$ for some $\s\in(k^{\times})^n$.
Then one can easily check that $E_j\cdot M_{(\s)}\subset
M_{(\epsilon_j \s)}$, where $1\leq j\leq n$,
$\epsilon_j=(q_j^{a_{j1}},q_j^{a_{j2}},\cdots,q_j^{a_{jn}})
\in(k^{\times})^n$ and
$\epsilon_j\s=(q_j^{a_{j1}}\s_1,q_j^{a_{j2}}\s_2,\cdots,q_j^{a_{jn}}\s_n)
$. It follows that $\sum_{\s\in(k^{\times})^n}M_{(\s)}$ is a
submodule of $M$, and the sum is a direct sum of vector spaces. Let
$\Pi (M)=\{\s \in (k^{\times})^n |M_{(\s)}\neq 0\}$, called the
weight space of $M$. If $M=\op _{\s \in(k^{\times})^n }M_{(\s)}$,
then $M$ is called a weight module.  Let $\mathcal {W}$ denote the
full subcategory of $_{\qgg}\mathcal {M}$ consisting of all the
weight modules. Obviously, $\mathcal {W}$ is closed under the direct
sum of modules.

Let $\s=(\s_1,\s_2,\cdots,\s_n)\in (k^{\times})^n$. Then $\s$
determines an algebra homomorphism $\s: \qgg\rightarrow k$ given by
$\s(K_i)=\s_i$ and $\s(E_i)=0$ for all $1\leq i\leq n$. Let $V_{\s}$
denote the corresponding $1$-dimensional $\qgg$-module. Then
$K_i\cdot v=\s_iv,\, E_i\cdot v=0$ for any $v\in V_{\s}$ and
$i=1,\cdots,n$. Obviously, $V_{\s}$ is a simple $\qgg$-module. For
any $\s$ and $\tau$ in $(k^{\times})^n$, $V_{\s}\cong V_{\tau}$ if
and only if $\s=\tau$.

Now let $M$ be a finite dimensional $\qgg$-module. Since $k$ is
algebraically closed, there is a $\s\in(k^{\times})^n$ such that
$M_{(\s)}\neq 0$. Hence $\op _{\s \in(k^{\times})^n }M_{(\s)}=\op
_{\s\in\Pi(M)}M_{(\s)}$ is a nonzero submodule of $M$. Thus if $M$
is a simple $\qgg$-module, then $M=\op _{\s\in\Pi(M)}M_{(\s)}$,
i.e., $M$ is a weight module. It follows that $M$ is a weight module
if $M$ is semisimple as a $\qgg$-module.

Now assume that $M$ is finite dimensional simple $\qgg$-module. Then
$\Pi(M)$ is a non-empty finite set. Since $q$ is not a root of
unity, there is a $\s\in\Pi(M)$ such that $\epsilon_j\s\notin\Pi(M)$
for all $1\leq j\leq n$. Hence $E_j\cdot M_{(\s)}=0$ for all $1\leq
j\leq n$, and any subspace of $M_{(\s)}$ is a submodule of $M$. It
follows that $M=M_{(\s)}\cong V_{\s}$. Thus we have the following
theorem.

\begin{theorem}\thlabel{2.1.1}
Any finite dimensional simple $\qgg$-module is a weight module and
must be $1$-dimensional. Moreover, there is a $1$-$1$
correspondences between $(k^{\times})^n$ and the set of isomorphism
classes of finite simple $\qgg$-modules.
\end{theorem}

\begin{remark}\relabel{2.1.2} For a Hopf algebra $H$, let $\mathcal{G}(H)$
denote the set of all the group-like elements in $H$. It is
well-known that an element $f\in H^*$ is a group-like element in
$H^\circ$ if and only if $f$ is an algebra map from $H$ to $k$ (cf.
\cite[Thm 9.1.4]{mon}). That is, $\mathcal {G}(H^\circ)={\rm
Alg}(H,k)$, the set of algebra homomorphisms from $H$ to $k$. By the
discussion above, an element $\s\in(k^{\times})^n$ determines an
algebra map $\s: \qgg\rightarrow k$, $\s(K_i)=\s_i$, $\s(E_i)=0$,
$1\leq i\leq n$. On the other hand, if $\s: \qgg\rightarrow k$ is an
algebra homomorphism, then one can easily check that $\s(K_i)\in
k^{\times}$ and $\s(E_i)=0$ for all $1\leq i\leq n$ since $q$ is not
a root of unity. Hence one can regard
$\s=(\s_1,\s_2,\cdots,\s_n)\in(k^{\times})^n$, where $\s_i=\s(K_i)$,
$i=1,2,\cdots,n$. With this identification, we have
$(k^{\times})^n={\rm Alg}(\qgg,k)=\mathcal {G}((\qgg)^\circ)$.
\end{remark}

\subsection{The Verma module}\selabel{2.2}
For any $\s \in (k^{\times})^n$, $V_{\s}$ is a $1$-dimensional
$\qg^0$-module since $\qg^0$ is a subalgebra of $\qgg$. Hence one
can form another $\qgg$-module $M(\s):=\qgg\ot_{\qg^0}V_{\s}$. We
call $M(\s)$ a Verma module. Since $\qgg\cong\qg^+\otimes\qg^0$ as
$\qg^+$-$\qg^0$-bimodules,
$\ms=\qgg\ot_{\qg^0}V_{\s}\cong(\qg^+\otimes\qg^0)\ot_{\qg^0}V_{\s}
\cong\qg^+\ot V_{\s}\cong\qg^+$ as $\qg^+ $-modules. That is, $\ms$
is a free $\qg^+$-module of rank $1$. Pick up a nonzero element
$x_{\s}\in V_{\s}$ and set $\vs=1\ot_{\qg^0}x_{\s}$ in $M(\s)$. Then
$\ms=\qg^+\cdot\vs$.

It is not difficult to check that $\ms \simeq M(\textbf{1})\ot
V_{\s}\simeq V_{\s}\ot M(\textbf{1})$ as left $\qgg$-modules for any
$\s \in (k^{\times})^n$, where $\textbf{1}=(1,1,\cds,1)\in
(k^{\times})^n$.

\begin{lemma}\lelabel{2.2.0} Let $\s$, $\tau \in (k^{\times})^n$. Then we have\\
$1)$ $\ms$ is a weight module;\\
$2)$ Let $J(\s)$ denote the submodule of $\ms$ generated by $E_1\cd
\vs,\cdots , E_n\cd \vs$. Then\\
\mbox{\hspace{0.4cm}}$J(\s)$ is a unique maximal submodule
of $\ms$;\\
$3)$ $\Ls:=\ms /J(\s)$ is a $1$-dimensional simple module and
$\Ls\cong V_{\s}$;\\
$4)$ $\ms$ is indecomposable.\\
$5)$ $\ms\cong M(\tau)$ if and only if $\s=\tau$.
\end{lemma}
\begin{proof}Note that $\qgg$ is a $Q_+$-graded Hopf algebra with
grading given by $deg(K_i)=0$ and $deg(E_i)=\a_i$, $i=1,2,\cdots,
n$. Obviously, $\qg^+$ is a graded subalgebra of $\qgg$. For any
$1\<i\<n$, the map $\qg^+\ra \qg^+$, $u\mapsto K_iuK_i^{-1}$, is a
graded algebra automorphism since $K_iE_jK_i^{-1}=q^{d_ja_{ji}}E_j$
for all $j=1,2,\cdots, n$. Now let $u\in\qg^+$ be a homogeneous
element with $deg(u)=\sum_{j=1}^nl_j\a_j\in Q_+=\Z_+\D$. Then
$K_iuK_i^{-1}=q^{\sum_{j=1}^nl_jd_ja_{ji}}u=q^{\sum_{j=1}^nl_jd_ia_{ij}}u$.
Hence in $M(\s)$, we have $K_i\cdot(u\cdot
v_{\s})=(K_iuK_i^{-1})\cdot(K_i\cdot
v_{\s})=\s_iq^{\sum_{j=1}^nl_jd_ia_{ij}}(u\cdot v_{\s})$ for all
$1\<i\<n$. Thus $u\cdot v_{\s}$ is a weight vector with the weight
$\tau\in(k^{\times})^n$ given by
$\tau_i=\s_iq^{\sum_{j=1}^nl_jd_ia_{ij}}$, $i=1,2,\cdots, n$. This
shows Part 1).

Let $u, u'\in\qg^+$. Since $C$ is a non-degenerate matrix (cf.
\cite{hum}) and $q$ is not a root of unity, it follows from the
above argument that $u\cdot v_{\s}\in M(\s)$ is a weight vector if
and only if $u$ is a homogeneous element. Moreover, if both $u$ and
$u'$ are homogeneous then $u\cdot v_{\s}$ and $u'\cdot v_{\s}$ have
the same weight if and only if $deg(u)=deg(u')$. Furthermore,
$u\cdot v_{\s}\in M(\s)_{(\s)}$ if and only if $deg(u)=0$ if and
only if $K_iuK_i^{-1}=u$ for all $1\<i\<n$.

From \cite[p.161]{jan}, we know that each $E_{\b_s}$ is a
homogeneous element with $deg(E_{\b_s})=\b_s\neq 0$, $s=1,
2,\cdots,N$. Hence if ${\bf k}=(k_1,k_2,\cdots,k_N)\in\Z_+^N$, then
$E^{\bf k}$ is a homogeneous element with $deg(E^{\bf
k})=\sum_{s=1}^Nk_s\b_s$. In particular, if ${\bf k}\neq 0$ then
$deg(E^{\bf k})\neq 0$. Since $\{E^{\bf k}|{\bf k}\in\Z_+^N\}$ is a
$k$-basis of $\qg^+$, the homogeneous component of $\qg^+$ of degree
0 is $k1$, where $1$ is the identity in $\qg^+$. It follows that
$M(\s)_{(\s)}=kv_{\s}$, $v_{\s}\notin J(\s)$ and
$\sum_{\tau\neq\s}M(\s)_{(\tau)}=\text{span}\{E^{\bf k}\cdot
v_{\s}|0\neq{\bf k}\in\Z_+^N\}$. Clearly, if ${\bf
k}=(k_1,k_2,\cdots,k_N)\neq 0$ in $\Z_+^N$, then $E^{\bf k}\cdot
v_{\s}\in J(\s)$. Therefore,
$J(\s)=\sum_{\tau\neq\s}M(\s)_{(\tau)}=\text{span}\{E^{\bf k}\cdot
v_{\s}|0\neq{\bf k}\in\Z_+^N\}$ and $M(\s)/J(\s)\cong V_{\s}$. This
shows Part 3).

It is easy to check that any submodule of a weight module is also a
weight module. Now let $M'$ be a proper submodule of $\ms$. Then
$M'$ is a weight module by Part (1). Since $v_{\s}\notin M'$ and
$M(\s)_{(\s)}=kv_{\s}$, $M'_{(\s)}=0$, i.e., $\s\notin\Pi(M')$.
Hence $M'=\sum_{\tau\in\Pi(M')}M'_{(\tau)}\subseteq J(\s)$ by the
last paragraph. This shows Parts 2) and 4).

If $\s=\tau$, then clearly $\ms\cong M(\tau)$. If $\s\neq\tau$, then
$V_{\s}\ncong V_{\tau}$, and so $\ms\ncong M(\tau)$. This shows Part
5).
\end{proof}

\begin{proposition}\prlabel{2.2.1}
Let $M$ be a $\qgg$-module and let $\s\in(k^{\times})^n$. Then there
is an epimorphism $f: \ms \ra M$ if and only if there exists a
weight vector $v\in M$ with weight $\s$ such that $M=\qgg \cdot
v=\qg^+\cdot v$.
\end{proposition}
\begin{proof}
It's trivial.
\end{proof}

Define a partial order "$\>$" on $(k^{\times})^n$: $\s \>\tau$ if
$\s\tau^{-1}=\epsilon_1^{s_1}\epsilon_2^{s_2}\cds \epsilon_n^{s_n}$
for some $s_1,s_2,\cds,s_n\in \Z_+$, where $\s
,\tau\in(k^{\times})^n$ and
$\s\tau^{-1}:=(\s_1\tau_1^{-1},\s_2\tau_2^{-1},\cdots,\s_n\tau_n^{-1})$.
Let $M$ be a $\qgg$-module. A weight vector $v\in M_{(\s)}$ with
weight $\s\in(k^{\times})^n$ is called a lowest weight vector if
there is no another weight vector $w$ with weight $\tau\in \Pi(M)$
such that $v\in \qgg \cd w$ and $w\notin \qgg \cd v$. A weight
module $M$ is called the lowest weight module if $M=\qgg \cd v$ for
some lowest weight vector $v$.

Clearly, every Verma module $M(\s)$ is a lowest weight module.
\begin{proposition}\prlabel{2.2.2}
Let $M$ is a weight module over $\qgg$. Then $M$ is a lowest weight
module if and only if $M$ is a quotient of some $\ms$.
\end{proposition}
\begin{proof}
It follows from \prref{2.2.1}.
\end{proof}

Note that $\mathcal {W}$ is the category of all the weight modules.
one can easily check that $\mcw$ is closed under tensor products,
submodules and quotient modules.

\begin{proposition}\prlabel{2.2.3}
Up to isomorphism of $\qgg$-modules, we have\\
$1)$ $\{V_{\s},\s \in(k^{\times})^n\}$ is a
complete set of simple objects in $\mcw$. \\
$2)$ $\{M(\s)|\s\in(k^{\times})^n\}$ is a set of nonisomorphic
indecomposable projective objects\\
\mbox{\hspace{0.5cm}}in $\mcw$.
\end{proposition}
\begin{proof}
1) Let $V$ be a simple object in $\mcw$. Since $V$ is a weight
module, one can pick up a nonzero weight vector $v\in V$ with weight
$\s$. Then $V=\qgg\cd v$, which is isomorphic to a quotient of $\ms$
by \prref{2.2.1}. It follows from \leref{2.2.0} that $V$ is
isomorphic to $V_{\s}$.

2) By \leref{2.2.0}, each $\ms$ is indecomposable. Let $f:M\ra L$ be
an epimorphism and $g:M(\s)\ra L$ be a morphism in $\mcw$. Then for
any $\tau\in(k^{\times})^n$, $f(M_{(\tau)})=L_{(\tau)}$ and
$g(M(\s)_{(\tau)})\subseteq L_{(\tau)}$. Hence there exists a weight
vector $m\in M_{(\s)}$ such that $f(m)=g(v_{\s})$. Define $\phi:
M(\s)\ra M$ by $\phi(u\cdot v_{\s})=u\cdot m$, $u\in \qg^+$. Since
$M(\s)$ is a free $\qg^+$-module with a $\qg^+$-basis $\{v_{\s}\}$,
$\phi$ is well-defined. It is easy to see that $\phi$ is a
$\qgg$-module morphism and $f\phi=g$. Hence $M(\s)$ is an
indecomposable projective object in $\mcw$ for any
$\s\in(k^{\times})^n$. By \leref{2.2.0}, $\ms\ncong M(\tau)$ if
$\s\neq\tau$. This completes the proof.
\end{proof}

Now we consider the tensor product of two weight modules. The
following lemma is obvious.

\begin{lemma}\lelabel{2.2.5}
Let $M$ and $N$ be two weight modules. Then $M\otimes N$ is also a
weight module and
$$(M\ot N)_{(\s)}=\bigoplus _{\tau \nu=\s}M_{(\tau)}\ot N_{(\nu)}.$$
\end{lemma}

Notice that $V_{\s}\ot V_{\tau}\cong V_{\tau}\ot V_{\s}\cong V_{\s
\tau}$ for any $\s, \tau \in (k^{\times})^n$.

We already know that $\qg^+$ is a $Q_+$-graded algebra, and the
homogeneous component $(\qg^+)_{\eta}$ of degree $\eta$ is equal to
$\mathrm{span}\{E_{\b_1}^{k_1}\cdots E_{\b_N}^{k_N}|\sum_{i=1}^N
k_i\b _i=\eta\}$, where $\eta \in Q_+$.

Define a group homomorphism $F:Q\ra (k^{\times})^n$ by
$F(\a_i)=\epsilon_i$, where $1\<i\<n$, and $\epsilon_i$ is given as
before. For a Verma module $\ms$, define a map $F_{\s}:Q_+\ra \Pi
(\ms)$ by
$$F_{\s}(\sum_{i=1}^nl_i\a_i)=F(\sum_{i=1}^nl_i\a_i)\s
=\epsilon_1^{l_1}\epsilon_2^{l_2}\cds\epsilon_n^{l_n}\s.$$ It is not
difficult to check that $F_{\s}$ is a bijective map. Let $\eta\in
Q_+$. If $F_{\s}(\eta)=\tau$, then
$\ms_{(\tau)}=(\qg^+)_{\eta}\cdot\vs$ by the proof of \leref{2.2.0},
and so ${\rm dim}\ms_{(\tau)}={\rm dim}(\qg^+)_{\eta}$.

For any $\s,\s ' \in (k^{\times})^n$,  $ \ms \ot M(\s ')$ is also a
weight module. We discuss the decomposition of $ \ms \ot M(\s ')$.

\begin{lemma}\lelabel{2.2.6}
Let $\s,\s ' \in (k^{\times})^n$. Then the elements $E^{\k}\cdot\vs
\ot v_{\s'}$ are lowest weight vectors in  $\ms \ot M(\s ')$, where
$\k \in \Z_+^N$.
\end{lemma}

\begin{proof} Let $\k \in \Z_+^N$ and $v=E^{\k}\cdot\vs \ot
v_{\s'}$. Then by the proof of \leref{2.2.0}, $v$ is a weight
vector. Assume there is a weight vector
$w=\sum_\textbf{r,s}a_{\textbf{r,s}} E^\textbf{r}\cdot\vs \ot
E^\textbf{s}\cdot v_{\s '}$ in $\ms \ot M(\s ')$ with weight $\tau$
such that $v\in \qgg \cd w$ and $w\notin \qgg \cd v$. Since $v$ is a
weight vector, $\qgg \cd v=\qg^+\cd v$. Hence there exists a
homogeneous element $h\in \qg^+$ such that
\begin{equation}
E^{\k}\cd\vs \ot v_{\s'}=h\cd(\sum_\textbf{r,s}a_{\textbf{r,s}}
E^\textbf{r}\cd\vs \ot E^{\textbf{s}}\cd v_{\s '})=\sum
_{\textbf{r,s}} a_{\textbf{r,s}} h_{(1)}E^\textbf{r}\cd\vs \ot
h_{(2)}E^\textbf{s}\cd v_{\s '}.\tag{$*$}
\end{equation}
Let $deg(h)=\eta\in Q_+$. Then $h\in(\qg^+)_{\eta}$. By
\cite[p.59]{jan}, we
have
$$\D((\qg^+)_{\eta})\subseteq\bigoplus_{0\<\mu\<\eta}
((\qg^+)_{\eta-\mu}K_{\mu}\otimes(\qg^+)_{\mu}.$$ From $w\notin \qgg
\cd v$, one knows that $\eta\neq 0$. Therefore, we have
$\D(h)=K_{\eta}\otimes h+\sum_{i=1}^mh_i\otimes g_i$, where
$h_i\in(\qg^+)_{\eta-\mu_i}K_{\mu_i}$ and $g_i\in(\qg^+)_{\mu_i}$
with $0\<\mu_i\lneq \eta$, $i=1,2,\cdots,m$. We may assume
$h_i=a_iE^{\k_i}K_{\mu_i}$ and $g_i=E^{{\bf l}_i}$ with $a_i\in
k^{\times}$ and $\k_i, {\bf l}_i\in\Z_+^N$ such that
$deg(E^{\k_i})=\eta-\mu_i$ and $deg(E^{{\bf l}_i})=\mu_i$,
$i=1,2,\cdots, m$. Thus from Eq.(*), one gets
$$E^{\k}\cd\vs \ot v_{\s'}=\sum_\textbf{r,s}a_{\textbf{r,s}}
K_{\eta}E^\textbf{r}\cd\vs \ot hE^{\textbf{s}}\cd v_{\s '}
+\sum_{\textbf{r,s},i} a_{\textbf{r,s}} h_iE^\textbf{r}\cd\vs \ot
g_iE^\textbf{s}\cd v_{\s '}.$$ By the proof of \leref{2.2.0}, we
have $a_{\textbf{r,s}} K_{\eta}E^\textbf{r}\cd\vs=a'_{\textbf{r,s}}
E^\textbf{r}\cd\vs$ with $a'_{\textbf{r,s}}=ba_{\textbf{r,s}}$ for
some $b\in k^{\times}$, and $a_{\textbf{r,s}}
h_iE^\textbf{r}\cd\vs=a_{\textbf{r,s}}a_iE^{\k_i}K_{\mu_i}E^\textbf{r}\cd\vs
=b_{\textbf{r,s},i}E^{\k_i}E^\textbf{r}\cd\vs$ with
$b_{\textbf{r,s},i}=c_ia_{\textbf{r,s}}$ for some $c_i\in
k^{\times}$. Therefore, we have
$$E^{\k}\cd\vs \ot v_{\s'}=\sum_\textbf{r,s}a'_{\textbf{r,s}}
E^\textbf{r}\cd\vs \ot hE^{\textbf{s}}\cd v_{\s '}
+\sum_{\textbf{r,s},i}b_{\textbf{r,s},i}E^{\k_i}E^\textbf{r}\cd\vs\ot
E^{{\bf l}_i}E^\textbf{s}\cd v_{\s '}.$$ Since $\ms$ and $M(\s')$
are two free $\qg^+$-modules with $\qg^+$-basis $v_{\s}$ and
$v_{\s'}$, respectively, $\ms\ot M(\s')$ is a free
$\qg^+\ot\qg^+$-module with a $\qg^+\ot\qg^+$-basis $v_{\s}\ot
v_{\s'}$. Hence the last equation is equivalent to
\begin{equation}
E^{\k}\ot 1=\sum_\textbf{r,s}a'_{\textbf{r,s}} E^\textbf{r}\ot
hE^{\textbf{s}}+
\sum_{\textbf{r,s},i}b_{\textbf{r,s},i}E^{\k_i}E^\textbf{r}\ot
E^{{\bf l}_i}E^\textbf{s}.\tag{$**$}
\end{equation}
Since $\{ E^{\k}|\k \in \Z_+^N\}$ is a $k$-basis of $\qg^+$ by
\thref{1.2.1}, $\{ E^{\k}\ot E^\textbf{l}|\k, \textbf{l}\in
\Z_+^N\}$ is a $k$-basis of $\qg^+\ot \qg^+$. It follows from
\cite[Corollary 1.8]{dck} that $\mathcal U$ has no zero divisors.
Hence $\{hE^{\textbf{s}}|\textbf{s}\in\Z_+^N\}$ are linearly
independent over $k$, and so are $\{E^\textbf{r}\ot
hE^{\textbf{s}}|\textbf{r, s}\in\Z_+^N\}$. On the other hand, since
$\qg^+$ is a $Q^+$-graded algebra, $\qg^+\ot\qg^+$ is a graded
$Q^+\oplus Q^+$-graded algebra with the grading given by $deg(
E^\textbf{r}\ot
E^{\textbf{s}})=(deg(E^\textbf{r}),deg(E^{\textbf{s}}))$, where
${\bf r, s}\in\Z_+^N$. When ${\bf s}=(s_1,s_2,\cdots,s_N)\neq{\bf
0}$ in $\Z_+^N$, $deg(E^{\bf s})=\sum_{j=1}^Ns_j\b_j\neq 0$. Let
$deg(E^{{\bf s}'})$ be a maximal element in the finite subset
$\{deg(E^{\bf s})|a'_{{\bf r,s}}\neq 0\text{ for some }{\bf r}\}$ of
$Q$. Then there is an ${\bf r}'$ such that $a'_{{\bf r}',{\bf
s}'}\neq 0$. We claim that $deg(E^{{\bf r}'}\ot hE^{{\bf
s}'})\notin\{deg(E^{\k_i}E^\textbf{r}\ot E^{{\bf
l}_i}E^\textbf{s})|b_{\textbf{r,s},i}\neq 0\}$. In fact, If
$deg(E^{{\bf r}'}\ot hE^{{\bf s}'})=deg(E^{\k_i}E^\textbf{r}\ot
E^{{\bf l}_i}E^\textbf{s})$ for some $b_{\textbf{r,s},i}\neq 0$,
then $deg(hE^{{\bf s}'})=deg(E^{{\bf l}_i}E^\textbf{s})$. Hence
$deg(h)+deg(E^{{\bf s}'})=deg(E^{{\bf l}_i})+deg(E^\textbf{s})$,
which implies that $deg(E^\textbf{s})-deg(E^{{\bf
s}'})=deg(h)-deg(E^{{\bf l}_i})=\eta-\mu_i\gneq 0$, a contradiction,
since $b_{\textbf{r,s},i}\neq 0$ iff $a_{{\bf r,s}}\neq 0$ iff
$a'_{{\bf r,s}}\neq 0$. Thus, the homogeneous component of degree
$deg(E^{{\bf r}'}\ot hE^{{\bf s}'})$ of the element
$\sum_\textbf{r,s}a'_{\textbf{r,s}} E^\textbf{r}\ot hE^{\textbf{s}}+
\sum_{\textbf{r,s},i}b_{\textbf{r,s},i}E^{\k_i}E^\textbf{r}\ot
E^{{\bf l}_i}E^\textbf{s}$ is not zero. Since $deg(hE^{{\bf
s}'})=\eta+deg(E^{{\bf s}'})\neq 0$, $deg(E^{{\bf r}'}\ot hE^{{\bf
s}'})\neq deg(E^{\k} \ot 1)=(\sum_{j=1}^Nk_j\b_j, 0)$. Thus, Eq.(**)
implies a contradiction. This completes the proof.
\end{proof}

\begin{theorem}\thlabel{2.2.7}
For any $\s,\s ' \in (k^{\times})^n$, there is a direct sum
decomposition in $\mcw$
\begin{equation}\eqlabel{5}
 \ms \ot M(\s ')
 =\bigoplus_{\k \in \Z_+^N} \qgg \cd(E^{\k}\cd\vs \ot v_{\s'}).
\end{equation}
\end{theorem}
\begin{proof}
Let $\k\neq{\bf l}\in \Z_+^N$. Assume $\qgg \cd(E^{\k}\cd\vs \ot
v_{\s'})\cap\qgg \cd(E^{\textbf{l}}\cd\vs \ot v_{\s'})\neq\{0\}$.
Then there is a nonzero weight vector $v$ in $\qgg \cd(E^{\k}\cd\vs
\ot v_{\s'})\cap\qgg \cd(E^{\textbf{l}}\cd\vs \ot v_{\s'})=\qg^+
\cd(E^{\k}\cd\vs \ot v_{\s'})\cap\qg^+ \cd(E^{\textbf{l}}\cd\vs \ot
v_{\s'})$. Hence $v=g\cd(E^{\k}\cd\vs \ot
v_{\s'})=h\cd(E^{\textbf{l}}\cd\vs \ot v_{\s'})$ for some
homogeneous elements $g,h\in \qg^+$. Let $deg(g)=\xi$ and
$deg(h)=\eta$ in $Q_+$. If $\eta=0$, then we may assume $h=1$, and
hence $E^{\textbf{l}}\cd\vs \ot v_{\s'}=g\cd(E^{\k}\cd\vs \ot
v_{\s'})\in\qg^+ \cd(E^{\k}\cd\vs \ot v_{\s'})$. From \leref{2.2.6},
both $E^{\k}\cd\vs \ot v_{\s'}$ and $E^{\textbf{l}}\cd\vs \ot
v_{\s'}$ are lowest weight vectors. Hence $E^{\k}\cd\vs \ot
v_{\s'}\in\qgg \cd(E^{\textbf{l}}\cd\vs \ot
v_{\s'})=\qg^+\cd(E^{\textbf{l}}\cd\vs \ot v_{\s'})$, and so
$E^{\k}\cd\vs \ot v_{\s'}=g'\cd(E^{\textbf{l}}\cd\vs \ot v_{\s'})$
for some homogeneous element $g'$ in $\qg^+$. Thus we have
$E^{\textbf{l}}\cd\vs \ot v_{\s'}=(gg')\cd(E^{\bf l}\cd\vs \ot
v_{\s'})$. By comparing their weights, one gets $deg(gg')=0$, which
implies $deg(g)=0$. So $g$ is a nonzero scale and
$E^{\textbf{l}}\cd\vs \ot v_{\s'}=g\cd(E^{\k}\cd\vs \ot
v_{\s'})=gE^{\k}\cd\vs \ot v_{\s'}$. Since $\ms \ot M(\s ')$ is a
free $\qg^+\oplus\qg^+$-module with a basis $\vs\otimes v_{\s'}$, it
follows that $E^{\textbf{l}}\ot 1=gE^{\k}\ot 1$, and so $\bf{l}=\k$,
a contradiction. Hence $\eta\neq 0$. Similarly, we also have
$\xi\neq 0$.

Now as in the proof of \leref{2.2.6}, we may assume
$$\begin{array}{rcl}
\D(h)&=&K_{\eta}\otimes h+\sum_{i=1}^ma_iE^{\k_i}K_{\eta_i}\otimes
E^{{\bf l}_i},\\
\D(g)&=&K_{\xi}\otimes g+\sum_{j=1}^nb_jE^{{\bf
s}_j}K_{\xi_j}\otimes
E^{{\bf r}_j},\\
\end{array}$$
where $a_i, b_j\in k^{\times}$, $0\<\eta_i\lneq \eta$,
$0\<\xi_j\lneq\xi$, $deg(E^{\k_i})=\eta-\eta_i$, $deg(E^{{\bf
l}_i})=\eta_i$, $deg(E^{{\bf s}_j})=\xi-\xi_j$ and $deg(E^{{\bf
r}_j})=\xi_j$. Hence we have
$$\begin{array}{rcl}
h\cd(E^{\textbf{l}}\cd\vs \ot v_{\s'})&=&aE^{\textbf{l}}\cd\vs \ot
h\cd v_{\s'}+\sum_{i=1}^ma_ic_iE^{\k_i}E^{\bf l}\cd\vs\ot E^{{\bf
l}_i}\cdot v_{\s'},\\
g\cd(E^{\k}\cd\vs \ot v_{\s'})&=&bE^{\k}\cdot\vs\ot g\cdot v_{\s'}
+\sum_{j=1}^nb_jd_jE^{{\bf s}_j}E^{\k}\cd\vs\ot E^{{\bf r}_j}\cd
v_{\s'}\\
\end{array}$$
for some $a, b, c_i, d_j\in k^{\times}$. Since $\ms \ot M(\s ')$ is
a free $\qg^+\oplus\qg^+$-module with a basis $\vs\otimes v_{\s'}$,
one gets
\begin{equation}
aE^{\textbf{l}}\ot h+\sum_{i=1}^ma_ic_iE^{\k_i}E^{\bf l}\ot E^{{\bf
l}_i}=bE^{\k}\ot g+\sum_{j=1}^nb_jd_jE^{{\bf s}_j}E^{\k}\ot E^{{\bf
r}_j}.\tag{$***$}
\end{equation}
in $\qg^+\ot\qg^+$. Note that $\qg^+\ot\qg^+$ is a $Q_+\oplus
Q_+$-graded algebra.  We have $deg(E^{\textbf{l}}\ot
h)=(deg(E^{\textbf{l}}), \eta)$, $deg(E^{\k_i}E^{\bf l}\ot E^{{\bf
l}_i})=(\eta-\eta_i+deg(E^{\bf l}), \eta_i)$, $deg(E^{\k}\ot
g)=(deg(E^{\k}),\xi)$ and $deg(E^{{\bf s}_j}E^{\k}\ot E^{{\bf
r}_j})=(\xi-\xi_j+deg(E^{\k}), \xi_j)$. Since $\eta\gneq \eta_i$ and
$\xi\gneq\xi_j$ for all $1\<i\<m$ and $1\<j\<n$, by considering the
degrees of the elements in Eq.(***), one gets $aE^{\textbf{l}}\ot
h=bE^{\k}\ot g$, which implies $E^{\textbf{l}}=E^{\k}$ and $h=cg$
for some $c\in k^{\times}$. This is impossible since ${\bf
l}\neq\k$. Thus we have proved that the sum $\sum_{\k \in \Z_+^N}
\qgg \cd(E^{\k}\cd\vs \ot v_{\s'})$ is direct.

From the above argument, one can see that $\qgg \cd(E^{\k}\cd\vs \ot
v_{\s'})=\qg^+\cd(E^{\k}\cd\vs \ot v_{\s'})$ is a free
$\qg^+$-module of rank 1 with a $\qg^+$-basis $E^{\k}\cd\vs \ot
v_{\s'}$, where $\k\in\Z_+^N$. Hence ${\rm
dim}((\qg^+)_{\eta}\cd(E^{\k}\cd\vs \ot v_{\s'}))={\rm
dim}(\qg^+)_{\eta}$ for any $\eta\in Q_+$ and $\k\in\Z_+^N$.

For any $\tau \in \Pi(\ms \ot M(\s '))$, the weight vector space
\begin{equation*}
\begin{split}
(\ms \ot M(\s '))_{(\tau)}
=&\bigoplus _{\nu \mu=\tau}\ms_{(\nu)}\ot M(\s')_{(\mu)}\\
=&\bigoplus_{\nu \mu=\tau}(\qg^+)_{F_{\s}^{-1}(\nu)}\cd\vs \ot
(\qg^+)_{F_{\s'}^{-1}(\mu)}\cd v_{\s '}\\
=&\bigoplus_{\substack{\eta,\xi\in
Q_+\\
F_{\s\s'}(\eta+\xi)=\tau}}(\qg^+)_{\eta}\cd\vs \ot
(\qg^+)_{\xi}\cd v_{\s '},\\
\end{split}
\end{equation*}
where the first equality follows from \leref{2.2.5}, the third
equality follows from the facts that $F_{\s}$ and $F_{\s'}$ are
bijective, and $F_{\s\s'}(\eta+\xi)=F_{\s}(\eta)F_{\s'}(\xi)$. On
the other hand, for any $\k=(k_1, k_2,\cdots,k_N),{\bf
l}=(l_1,l_2,\cdots,l_N) \in \Z_+^N$, $E^{\bf l}\cd(E^{\k}\vs \ot
v_{\s'})$ is a weight vector with weight $F_{\s\s'}(deg(E^{\bf
l})+deg(E^{\k}))=F_{\s\s'}(deg(E^{\bf
l}E^{\k}))=F_{\s\s'}(\sum_{j=1}^N(l_j+k_j)\b_j)$. Hence we have
\begin{equation*}
\begin{split}
&(\bigoplus_{\k \in \Z_+^N} \qgg \cd(E^{\k}\cd\vs \ot
v_{\s'}))_{(\tau)}\\
=&(\bigoplus_{\k \in \Z_+^N}
\qg^+ \cd(E^{\k}\cd\vs \ot v_{\s'}))_{(\tau)}\\
=&(\bigoplus_{\k,{\bf l} \in \Z_+^N}
E^{\bf l}\cd(E^{\k}\cd\vs \ot v_{\s'}))_{(\tau)}\\
=&\text{span}\{E^{\textbf{l}}\cd (E^{\k}\cd\vs \ot v_{\s'})|\k,{\bf
l}
\in \Z_+^N, F_{\s\s'}(deg(E^{\bf l})+deg(E^{\k}))=\tau\}\\
=&\bigoplus_{\substack{\eta,\xi\in
Q_+\\
F_{\s\s'}(\eta+\xi)=\tau}}(\qg^+)_{\eta}\cd((\qg^+)_{\xi}\cd\vs \ot
v_{\s '}) .\\
\end{split}
\end{equation*}

Now let $\eta, \xi\in Q_+$ with $F_{\s\s'}(\eta+\xi)=\tau$. Since
$\ms \ot M(\s ')$ is a free $\qg^+\oplus\qg^+$-module with a basis
$\vs\otimes v_{\s'}$, one knows that ${\rm dim}((\qg^+)_{\eta}\cd\vs
\ot (\qg^+)_{\xi}\cd v_{\s '})={\rm dim}((\qg^+)_{\eta}\ot
(\qg^+)_{\xi})=({\rm dim}(\qg^+)_{\eta})({\rm
dim}(\qg^+)_{\xi})<\infty$. On the other hand, since
$\{E^{\k}|\k\in\Z_+^N, deg(E^{\k})=\xi\}$ is a $k$-basis of
$(\qg^+)_{\xi}$, we have
\begin{equation*}
\begin{split}
(\qg^+)_{\eta}\cd((\qg^+)_{\xi}\cd\vs \ot v_{\s '})
&=\sum_{\substack{\k\in\Z_+^N\\
deg(E^{\k})=\xi}}(\qg^+)_{\eta}\cd(E^{\k}\cd\vs \ot v_{\s '})\\
&=\bigoplus_{\substack{\k\in\Z_+^N\\
deg(E^{\k})=\xi}}(\qg^+)_{\eta}\cd(E^{\k}\cd\vs \ot v_{\s '}),\\
\end{split}
\end{equation*}
and hence,
\begin{equation*}
\begin{split}
{\rm dim}((\qg^+)_{\eta}\cd((\qg^+)_{\xi}\cd\vs \ot
v_{\s'}))&=\sum_{\substack{\k\in\Z_+^N\\
deg(E^{\k})=\xi}}{\rm dim}((\qg^+)_{\eta}\cd(E^{\k}\cd\vs\ot v_{\s'}))\\
&=({\rm dim}(\qg^+)_{\xi})({\rm dim}(\qg^+)_{\eta}).\\
\end{split}
\end{equation*}
Thus, ${\rm dim}((\qg^+)_{\eta}\cd\vs \ot (\qg^+)_{\xi}\cd v_{\s
'})={\rm dim}((\qg^+)_{\eta}\cd((\qg^+)_{\xi}\cd\vs \ot v_{\s'}))$.
Since $\{(\eta, \xi)|\eta, \xi\in Q_+, F_{\s\s'}(\eta+\xi)=\tau\}$
is a finite set, we have ${\rm dim}((\ms \ot M(\s '))_{(\tau)})={\rm
dim}((\bigoplus_{\k \in \Z_+^N} \qgg \cd(E^{\k}\vs \ot
v_{\s'}))_{(\tau)})$, and so $(\ms \ot M(\s
'))_{(\tau)}=(\bigoplus_{\k \in \Z_+^N} \qgg \cd(E^{\k}\vs \ot
v_{\s'}))_{(\tau)}$. This completes the proof of the theorem.
\end{proof}

\begin{corollary}\colabel{}
For any $\s,\s ' \in (k^{\times})^n$, there is a $\qgg$-module
isomorphism
$$ \ms \ot M(\s ')\simeq \bigoplus _{\eta\in Q^+}({\rm dim}(\qg^+)_{\eta})M(F(\eta)\s\s').$$
\end{corollary}
\begin{proof}
For $\k =(k_1,k_2,\cds,k_N)\in \Z_+^N$, let $\eta=deg(E^{\k})=\sum
_{j=1}^Nk_j\b_j=\sum_{i=1}^n s_i\a_i$ in $Q_+$. Then from the proof
of \thref{2.2.7}, one gets that $\qgg \cd(E^{\k}\vs \ot v_{\s'}) $
is isomorphic to $M(F(\eta)\s\s')$, where
$F(\eta)\s\s'=F_{\s\s'}(\eta) =\epsilon_1^{s_1}\epsilon_2^{s_2}\cds
\epsilon_n^{s_n}\s\s'$. The corollary follows from \thref{2.2.7}.
\end{proof}

\begin{corollary}\colabel{}
$ \ms \ot M(\s ')\simeq M(\s ')\ot \ms$ as left $\qgg$-modules for
all $\s,\s ' \in (k^{\times})^n$.
\end{corollary}

\subsection{The comodule and the Yetter-Drinfel'd module}\selabel{2.3}
Recall that the map $ht:Q\ra \Z$, $ht(\sum_{i=1}^na_i\a
_i)=\sum_{i=1}^na_i$, is $\Z$-linear. That is, $ht$ is a group
homomorphism from the abelian group $Q$ to the abelian $\Z$. Thus
one can define a $\Z_+$-grading $h$ on $\qgg$ by
$$h(E^{\k}K_{\l})=\sum _{j=1}^N k_j ht(\b_j),$$
where $E^{\k}K_{\l}=E_{\b_1}^{k_1}\cdots E_{\b_N}^{k_N}
K_1^{t_1}\cdots K_n^{t_n}$, ${\bf k}=(k_1,\,\cdots,\,k_N)\in \Z_+^N$
and $\l=\sum_{i=1}^nt_i\a_i\in Q$. Then $\qgg
=\bigoplus_{l=0}^{\infty}(\qgg)_{(l)}$ is a $\Z _+$-graded Hopf
algebra, where
$(\qgg)_{(l)}=\mathrm{span}\{E^{\k}K_{\l}|h(E^{\k}K_{\l})=l,\k
\in\Z_+^N,\l\in Q \}$. Particularly, $(\qgg)_{(0)}=k[ K_1^{\pm
1},K_2^{\pm },\cds,K_n^{\pm}]$ is a Laurent polynomial algebra in
$n$ variables.

Since the coradical of $\qgg$ is contained in $(\qgg)_{(0)}$(cf.
\cite[Lemma 5.3.4]{mon}), one can see that
$(\qgg)_{(0)}=k\mathcal{G}(\qgg)$ is the coradical of $\qgg$, and
$\mathcal{G}(\qgg)=\{K_{\l}|\l\in Q \}$. Hence each simple right
coideal of $\qgg$ is a simple subcoalgebra and has the form $kg$,
where $g\in \mathcal{G}(\qgg)$.

Let $H$ be a Hopf algebra with a bijective antipode $S$. Recall that
a (left-right) Yetter-Drinfel'd $H$-module (simply, YD $H$-module)
$M$ is a triple $(M,\cd,\rho )$ such that $(M, \cd)$ is a left
$H$-module, $(M,\rho)$ is a right $H$-comodule, and the following
equivalent compatibility conditions are satisfied (cf. \cite{lr}):
\begin{eqnarray}
h_{(1)}\cd m_{(0)}\ot h_{(2)}m_{(1)} =(h_{(2)}\cd m)_{(0)}\ot
(h_{(2)}\cd m)_{(1)}h_{(1)},
\end{eqnarray}
\begin{eqnarray}\eqlabel{15}
\rho (h\cd m)=h_{(2)}\cd m_{(0)}\ot h_{(3)}m_{(1)}S^{-1}(h_{(1)}),
\end{eqnarray}
where $h\in H$ and $m\in M$. A YD $H$-module map between two YD
$H$-modules is simultaneously an $H$-module map and an $H$-comodule
map. Let $_H  \yd ^H$ denote the category of (left-right) YD
$H$-modules and YD $H$-module maps. This is a braided monoidal
category. Let $M\in\ _H\!\yd ^H$. A YD $H$-submodule of $M$ is both
an $H$-submodule and an $H$-subcomodule of $M$. A YD $H$-module is
simple if it has no non-trivial YD $H$-submodules.

In \cite{rad}, Radford provided a procedure to construct some YD
$H$-modules through modules and comodules over $H$.

Let $L$ be a simple module over $H$. Then $L\ot H \in _H\!\! \yd ^H$
with the action and coaction of $H$ given by
$$h\cd (l\ot a)=h_{(2)}\cd l\ot h_{(3)}aS^{-1}(h_{(1)}),\quad \rho (l\ot
a)=(l\ot a_{(1)})\ot a_{(2)},\,$$ where $l\in L,h,a\in H$ (see
\cite{rad}). Clearly, $L\ot H \simeq(\mathrm{dim}L) H$ as comodules
over $H$.

For $\b \in \mathcal {G}(H^\circ)$, let $L_\b$ denote the
corresponding $1$-dimensional simple $H$-module: $h\cdot l=\b(h)l$,
$h\in H$, $l\in L_{\b}$. Define $H_\b =(H, \cd _\b ,\D)\in _H\!\!
\yd ^H$, where $h\cd _\b a=(h_{(2)}\!\!\leftharpoonup \b)
aS^{-1}(h_{(1)})=\b(h_{(2)})h_{(3)}aS^{-1}(h_{(1)})$. Then $H_\b
\simeq L_\b \ot H$. Let $N$ be a right coideal of $H$. Then
$H_{\b,N}=H\cd _\b N$ is a YD $H$-submodule of $H_\b$ (see
\cite{rad}).

\begin{lemma}\lelabel{2.3.1}$($\cite{rad}$)$
Let $H=\op _{n=0}^ \infty H_{(n)}$ be a graded bialgebra over $k$
which affords $H^{op}$ the structure of a graded Hopf algebra over
$k$. Suppose that $H_{(0)}$ is a commutative cocommutative Hopf
subalgebra of $H^{op}$.
\begin{itemize}
\item[a)] Let $\b: H\ra k$ be a map of graded algebras and let $N$ be a
simple right coideal of $H$. Then $H_{\b ,N}$ is a simple YD
$H$-module.
\item[b)] Suppose $\b, \b ':H\ra k $ are maps of graded algebras
and $ N,N'$ are simple coideals of $H$. Then $H_{\b ,N}\simeq H_{\b
' ,N'}$ if and only if $\b=\b '$ and $N=N'$.
\end{itemize}
\end{lemma}

We will apply these results to the half quantum group $\qgg$.

\begin{lemma}\lelabel{}
Let $\b: \qgg \ra k$ be an algebra map. Then $\b$ is a map of graded
algebras.
\end{lemma}

\begin{proof} It follows from \reref{2.1.2}.
\end{proof}

For any $\b \in \mcg ((\qgg)^\circ)$ and $g\in \mcg (\qgg)$, we have
$(\qgg)_{\b,g}=\qgg\cd _\b g\subset (\qgg)_{\b}$, where $h\cd _\b
g=\b(h_{(2)})h_{(3)}gS^{-1}(h_{(1)})$, $h\in \qgg$ and $S$ is the
bijective antipode of $\qgg$. By \leref{2.3.1}, $(\qgg)_{\b,g}$ is a
simple YD $\qgg$-module since $kg$ is a simple right coideal of
$\qgg$.

\begin{theorem}\thlabel{}
Let $\Phi(\b,g)=[(\qgg)_{\b,g}]$. Then $\Phi$ is a bijective map
from $\mcg ((\qgg)^\circ)\times \mcg (\qgg)$ to the set
$\mathcal{E}$ of isomorphic classes  of simple {\rm YD}
$\qgg$-modules which are weight modules as $\qgg$-modules.
\end{theorem}
\begin{proof}
From the discussion before, $(\qgg)_{\b,g}$ is a simple YD
$\qgg$-module. It is easy to see that $(\qgg)_{\b,g}$ is also a
weight $\qgg$-module. By \leref{2.3.1}, $\Phi$ is injective.

Now let $M$ be a simple YD $\qgg$-module which is a weight module.
Then $M=\op _{\s\in(k^{\times})^n}M_{(\s)}$. Using an argument
similar to \cite[p.697]{rad}, one can show that $M$ is a
Yetter-Drinfel'd $\qgg$-submodule of some $L\ot \qgg$. In fact, let
$N$ be a simple $\qgg$-subcomodule of $M$. Then $N=km$ with $\rho
(m)=m\ot g$ for some $m\in M$ and $g\in \mcg (\qgg)$ since each
simple subcoalgebra of $\qgg$ has the form $kg$, $g\in \mcg (\qgg)$.
Let $m=\sum _{j=1}^sm_j$ with $0\neq m_j\in M_{(\tau ^j)}$, where
$\tau ^1,\tau^2,\cds,\tau^s$ are distinct elements in $\Pi(M)$. If
$s\>2$, then $\tau_i^1\neq \tau_i^2$ for some $1\<i\<n$. By
Eq.\equref{15}, $\rho (K_i\cd m)=K_i\cd m\ot K_i gS(K_i)=K_i\cd m\ot
g$, i.e., ${\rm span}\{K_i\cd m\}$ is also a simple subcomodule
which is isomorphic to $N$. Let $m'=\tau_i^1m-K_i \cd m$. Then $km'$
is a simple subcomodule of $M$ and $m'\in\bigoplus _{j=2}^s M_{(\tau
^j)}$. Continuing this process, one may assume that $m$ is a weight
vector with weight $\s$. Note that $\qgg \cd N$ is a YD
$\qgg$-submodule of $M$ by Eq.\equref{15}. So $M=\qgg \cd N$ and $M$
is isomorphic, as a $\qgg$-module, to some quotient of $\ms$ by
\prref{2.2.1}. Since $M(\s)$ has a unique maximal $\qgg$-submodule,
so does $M$. Let $M'$ be the unique maximal $\qgg$-submodule of $M$.
Then $L=M/M'$ is a $1$-dimensional $\qgg$-module and $L\cong
V_{\s}$. Define $f:M\ra L\ot \qgg$ by $f(m)=p(m_{(0)})\ot m_{(1)}$,
where $p:M\ra L$ denotes the canonical epimorphism.  $L\ot\qgg$ is a
YD $\qgg$-module with the $\qgg$-action and $\qgg$-coaction given as
before. It is easy to see that $f$ is a YD $\qgg$-module map. Since
$M$ is a simple YD $\qgg$-module and $f\neq 0$, $\mathrm{Ker}f=0$,
and so $M$ is a YD $\qgg$-submodule of $L\ot \qgg$.

Note that $L\ot\qgg\cong V_{\s}\ot\qgg\cong (\qgg)_\b$, where
$\b=\s\in(k^{\times})^n={\rm Alg}(\qgg, k)=\mcg((\qgg)^{\circ})$ as
stated in \reref{2.1.2}. Hence we may regard $M\subset (\qgg)_\b$.
Thus $N=kg$ for some $g\in \mcg (\qgg)$, and so $M=\qgg \cd _\b
g=(\qgg)_{\b, g}$. It follows that $\Phi$ is surjective.

\end{proof}

\end{document}